\documentclass[11pt,a4paper]{article}
\usepackage[T1]{fontenc}
\usepackage{lmodern}
\usepackage{amsmath,amssymb,amsthm}
\usepackage[margin=29mm]{geometry}
\usepackage[hidelinks]{hyperref}
\hypersetup{
  pdftitle={Negative Ricci Curvature under Direct Sums and Levi Extensions},
  pdfsubject={Left-invariant Riemannian metrics on real Lie groups},
  pdfkeywords={negative Ricci curvature, Levi decomposition, nilradical, Lie groups}
}
\allowdisplaybreaks[2]
\numberwithin{equation}{section}

\newtheorem{theorem}{Theorem}[section]
\newtheorem{proposition}[theorem]{Proposition}
\newtheorem{lemma}[theorem]{Lemma}
\newtheorem{corollary}[theorem]{Corollary}
\theoremstyle{definition}

\newtheorem{example}[theorem]{Example}
\theoremstyle{remark}
\newtheorem{remark}[theorem]{Remark}

\newcommand{\Ric}{\operatorname{Ric}}
\newcommand{\tr}{\operatorname{tr}}

\newcommand{\End}{\operatorname{End}}
\newcommand{\Sym}{\operatorname{Sym}}
\newcommand{\ad}{\operatorname{ad}}
\newcommand{\RR}{\mathbb R}
\newcommand{\CC}{\mathbb C}
\newcommand{\uu}{\mathfrak u}
\newcommand{\lln}{\mathfrak l_-}
\newcommand{\ssplit}{\mathfrak s}

\newcommand{\bb}{\mathfrak b}
\newcommand{\aaa}{\mathfrak a}

\newcommand{\ttc}{\mathfrak t}
\newcommand{\Sop}{\operatorname{S}}

\title{Negative Ricci Curvature under Direct Sums\\ and Levi Extensions}
\author{}
\date{}
\begin{document}
\maketitle

\begin{abstract}
We prove that adjoining any finite number of copies of
\(\mathfrak{sl}(2,\mathbb R)\) to a nonzero real Lie algebra with a
negative Ricci inner product preserves the existence of such an inner
product. The construction preserves the prescribed metric on the
original summand. In particular,
\(\mathrm{SL}(2,\mathbb R)\times\mathrm{SL}(3,\mathbb R)\) has a
left-invariant metric of strictly negative Ricci curvature, although
its three-dimensional factor has none. We use the construction to
prove a criterion for Lie algebras with a stable radical direction
and a Levi algebra consisting of a compact part, a part already
carrying a negative Ricci inner product, and split three-dimensional
simple factors. In this class the criterion is faithfulness of the
compact action on the actual nilradical. We determine the sharp
dimension and the equality structures for a prescribed Levi algebra.
We also prove a converse to the direct-sum construction for three
types of solvable summand and classify the nonsolvable, nonsemisimple
examples in dimensions five and six.
\end{abstract}

\section{Introduction}\label{sec:introduction}

Can a semisimple Lie group carry a left-invariant Riemannian
metric of strictly negative Ricci curvature if one of its
simple factors carries no such metric? Lauret and Will asked
the corresponding factorwise existence question
\cite[Section~2.2]{LWSurvey}. The implication from the simple
factors to the product follows by taking a Riemannian
product. In the reverse direction, a left-invariant metric
on the product may pair vectors from distinct simple ideals.
We use such pairings to answer the question and then use the
same construction to treat factors acting trivially in a
Levi extension.

Throughout, Lie algebras are real and finite dimensional.
A negative Ricci inner product is a positive definite inner
product whose full Levi--Civita Ricci tensor is strictly
negative definite.

\subsection*{Direct sums and the semisimple factor question}

Leite and Dotti de Miatello constructed negative Ricci metrics
on \(\mathrm{SL}(n,\RR)\) for \(n\geq3\) \cite{DL}.
Dotti Miatello, Leite and Miatello subsequently gave
constructions for broad classes of complex simple and split
real simple groups \cite{DLM}. Milnor's three-dimensional
curvature formulas exclude a negative Ricci metric on
\(\mathfrak{sl}(2,\RR)\) itself
\cite[Theorem~4.3 and Corollary~4.7]{Milnor}. A compact
semisimple ideal gives a different obstruction, valid for
arbitrary pairings between ideals
\cite[Lemma~1.6]{JP}, \cite[Proposition~2.1]{AraujoObstruction}.
The split three-dimensional simple algebra thus provides
a test of factorwise inheritance among noncompact factors.

Our direct-sum theorem (Theorem~\ref{thm:attachment}) starts
with any nonzero Lie algebra \(\mathfrak h\) carrying a
prescribed negative Ricci inner product. It gives a negative
Ricci inner product on
\[
 \mathfrak{sl}(2,\RR)^{\oplus m}\oplus\mathfrak h
 \qquad (m\geq1)
\]
whose restriction to \(\mathfrak h\) is the prescribed
inner product. For one added factor, we choose a unit vector
\(v\in\mathfrak h\) with \(\tr\ad_v=0\) and add a small
multiple \(tv\) to one vector in a frame for the split
factor. The change gives an independent Ricci eigenvalue
\(-3t^2/8\) in the formerly obstructing direction. The
remaining mixed entries are of order \(t\), and their
contribution to a scalar Schur complement is of order
\(t^2\). This allows the negative term of order \(t\)
in that complement to control the sign.

Taking \(\mathfrak h=\mathfrak{sl}(3,\RR)\) yields a negative
Ricci metric on
\(\mathrm{SL}(2,\RR)\times\mathrm{SL}(3,\RR)\)
(Corollary~\ref{cor:attachment-semisimple}). The
eleven-dimensional example answers the semisimple factor
question in the negative. The same argument works with
\(\mathfrak{sl}(n,\RR)\) for every \(n\geq3\). Recent results
about cyclic metrics concern a class in which split
three-dimensional factors split off as Riemannian products
\cite[Theorem~3.4]{Cyclic}. The metric in our construction
has \(\langle[k,q],v\rangle=t\) in the frame of
\eqref{eq:attachment-brackets}, with \([q,v]=[v,k]=0\),
so the cyclic sum is nonzero. For comparison, the
isometric splitting theorem of B\"ohm and Lafuente concerns
negative Einstein metrics. Their result for negative Ricci
curvature along the orbits of a unimodular group gives
semisimplicity of that group
\cite[Theorems~A and~4.1]{BLUnimodular}.
The \(\mathrm{SL}(3,\RR)\) metrics of \cite{DL} have
also been used to study conjugate points along a geodesic
with negative directional Ricci curvature
\cite[Example~4.12, arXiv v4]{LLP}.

There is an earlier five-dimensional instance of the
direct-sum phenomenon. Andrada, Fino and Vezzoni identify the
algebra in their classification of Sasakian Lie algebras
\cite{AFV}. Boeckx's five-dimensional Lie group models
with invariant \(I<-1\) have this algebra, as noted
by Loiudice and Lotta
\cite{Boeckx}, \cite[Remark~1]{LL}. The contact family
constructed by Loiudice and
Lotta on \(\mathrm{SL}(2,\RR)\times\mathrm{Aff}(\RR)\),
together with the Ricci formulas for contact
\((\kappa,\mu)\)-spaces, gives negative Ricci metrics for an
interval of parameters \cite[Theorem~6]{LL},
\cite[Equations~(2.6)--(2.7)]{LottaMM}. We give the exact
interval for that family in Section~\ref{sec:lowdim}.
The direct-sum theorem extends the existence phenomenon to
an arbitrary negative Ricci summand and applies when that
summand is semisimple.

\subsection*{Stable directions and Levi actions}

Write \(\mathfrak g=\mathfrak b\ltimes\mathfrak r\) for a
Levi decomposition and \(\mathfrak n\) for the actual
nilradical. An element \(Y\in\mathfrak r\) is called a
\emph{stable radical direction} if every eigenvalue of
\(\ad_Y|_{\mathfrak n}\) has strictly positive real part.
For solvable Lie algebras, Nikolayevsky and Nikonorov proved
that the existence of such a direction suffices for
negative Ricci curvature \cite[Theorem~2]{NN}. The
rank-one negative sectional curvature construction goes
back to Heintze \cite[Theorem~3]{Heintze}.
For rank-one solvable extensions, Guti\'errez studies the
cone of Ricci-negative derivations and its relation to the
moment map \cite{Gutierrez}.

For semisimple actions, a related result announced by
Dotti Miatello treats an abelian radical and gives a
nonpositive Ricci metric from a nonpositive Ricci metric
on the noncompact semisimple factor adapted to a Cartan
decomposition \cite[p.~223]{DottiAnnouncement}. See also
her subsequent paper on nonpositive Ricci metrics on
semidirect products \cite{DottiSemidirect}. In the strict
case, Will works with an abelian complement commuting with
the semisimple algebra, a common normal action of that
complement, and a stable direction. Her hypotheses also
require that no nonzero vector in the complement have
entirely imaginary spectrum and that the semisimple algebra
carry a negative Ricci metric making a Cartan decomposition
orthogonal \cite[Theorem~5.2, arXiv v2]{WillII}.

For compact Levi factors, Will's first paper treats
scalar extensions of nontrivial real
\(\mathfrak{su}(2)\)-modules \cite{WillI}. Emilio Lauret and
Will obtain scalar-extension metrics from an explicit
condition on the weight spaces of a real module. In
particular, they cover all but finitely many complex
irreducible modules for each compact simple algebra, regarded
as real, and extend the method to certain actions on
nonabelian nilradicals
\cite[Theorems~1.1, 3.1 and~5.3]{LW}. Their Remark~4.3
identifies small complex modules left out by their three
weight constructions. These results motivate a criterion
in terms of the action on \(\mathfrak n\), allowing an
arbitrary faithful real compact action and several radical
directions.

\begin{theorem}[Stable radical criterion]\label{thm:mixed}
Let \(\mathfrak g=\mathfrak b\ltimes\mathfrak r\) have
Levi algebra
\[
 \mathfrak b=\mathfrak u\oplus\mathfrak l_-
                \oplus\mathfrak{sl}(2,\RR)^{\oplus k},
 \qquad k\geq0 .
\]
Suppose \(\mathfrak u\) is compact semisimple and the
semisimple algebra \(\mathfrak l_-\) admits a negative
Ricci inner product. Either summand may be zero.
Suppose that \(\mathfrak n\ne0\) and that
\(\mathfrak r\) contains a stable radical direction.
Then \(\mathfrak g\) admits a negative Ricci inner product
if and only if the action of \(\mathfrak u\) on
\(\mathfrak n\) is faithful.
\end{theorem}

The necessity follows from the compact-ideal obstruction.
For sufficiency, equivariant bracket limits first reduce
the problem to an abelian nilradical with a joint character
decomposition. A deformation along compact Cartan weight
planes then controls the compact Ricci block and its pairing
with the commuting semisimple action. A second deformation
controls the nontrivially acting split factors. The
complement can have directions with entirely imaginary
characters. A small real perturbation gives a definite
Ricci block in those directions and the trace identities
in \eqref{eq:anti-cross} control the pairing with the Levi
algebra. The direct-sum theorem adjoins each split factor whose
action on \(\mathfrak n\) is trivial.

In particular, the radical may have a complement of dimension
greater than one, nonreal spectrum or nontrivial Jordan
blocks. If the Levi algebra itself already carries a negative
Ricci metric, the proof can start with the prescribed metric
(Corollary~\ref{thm:negative-base}). The proof of the mixed
criterion uses separate features of the three types of
Levi ideal, and requires negative Ricci curvature on all
of \(\mathfrak l_-\) as a hypothesis.

\subsection*{Dimension bounds, converse results and explicit metrics}

For a fixed Levi algebra \(\mathfrak b\) as in the theorem
with \(\mathfrak u\ne0\), the minimum dimension is
\(\dim\mathfrak b+1+\mu_{\RR}(\mathfrak u)\), where
\(\mu_{\RR}(\mathfrak u)\) is the least dimension of a
faithful real \(\mathfrak u\)-module. The lower bound holds
for arbitrary radicals, and we determine the equality
structures in Theorem~\ref{thm:minimum-dimensions}. For
types \(A_1,A_2,B_2,G_2\) with compact simple Levi algebra
the minima are \(7,15,16,22\). The seven-dimensional compact
\(A_1\) example belongs to the earlier constructions
\cite{WillI}. Table~1 of \cite{LWSurvey} concerns a fixed
construction using complex irreducible modules viewed as
real, so its dimension convention differs from minimization
over all faithful real modules. Proposition~\ref{prop:scalar-construction}
gives one metric formula and an explicit negative upper
bound for the Ricci tensor in the scalar case. Anisotropic
compact constructions already appear in
\cite[Example~6.11 and Remark~6.12]{LWSCF}. Araujo's
Ricci block formulas also give the scalar-case existence
by taking an anisotropic metric on the module and then
rescaling the dilation direction
\cite[Proposition~4.6 and Corollary~4.8]{AraujoFlow}.

We also prove reverse inheritance from
\(\mathfrak{sl}(2,\RR)^{\oplus m}\oplus\mathfrak r\) to the
solvable factor \(\mathfrak r\) when its actual nilradical is
abelian, Heisenberg, or standard filiform. Tracing the full
Ricci tensor on abelian ideals recovers the already known
solvable-side criteria \cite[Theorems~2--4]{NN}.
Nikolayevsky's later work treats the solvable problem for
general filiform nilradicals \cite{NikFiliform}. These
conditions and a direct calculation give the five- and
six-dimensional classifications in
Theorem~\ref{thm:lowdim}. Turkowski's classification provides
the underlying low-dimensional algebra structures
\cite[Table~I]{Turkowski}. In dimension six we single out
the natural real affine algebra
\(\mathfrak{gl}(2,\RR)\ltimes\RR^2\). Will's related
construction in \cite[Remark~4.9, arXiv v2]{WillII} uses the
complex standard module regarded as real, which gives an
eight-dimensional algebra. The explicit metric here uses
the real two-dimensional standard module. Smolentsev and
Sokolova compute curvature on this same affine algebra
for indefinite K\"ahler and semi-para-K\"ahler structures
\cite[Section~3.1]{SmolentsevSokolova}.

\section{Ricci formulas and structural obstructions}\label{sec:prelim}

For a metric Lie algebra \(\mathfrak g\), write \(B_{\mathfrak g}\)
for its Killing form, \(A^*\) for the metric adjoint of an
endomorphism \(A\), and \(\Sop(A)=(A+A^*)/2\). The mean-curvature
vector \(H\) satisfies \(\langle H,x\rangle=\tr\ad_x\).
For an orthonormal basis \((e_i)\), the Ricci form is
\begin{align}
 \Ric(x,y)={}&-\frac12\sum_i\langle[x,e_i],[y,e_i]\rangle
 +\frac14\sum_{i,j}\langle[e_i,e_j],x\rangle
                       \langle[e_i,e_j],y\rangle \notag\\
 &-\frac12 B_{\mathfrak g}(x,y)
   -\langle\Sop(\ad_H)x,y\rangle .
 \label{eq:general-ricci}
\end{align}
We use the same symbol for the corresponding self-adjoint operator
when writing operator formulas. See \cite{LWSurvey,WillII}.

\begin{lemma}\label{lem:radical}
Let \(\mathfrak g=\mathfrak b\ltimes\mathfrak r\) be a Levi
decomposition and let \(\mathfrak n\) be the actual nilradical.
Then \(\mathfrak n=\operatorname{nilrad}(\mathfrak r)\) and
\([\mathfrak g,\mathfrak r]\subseteq\mathfrak n\). If an ideal
\(\mathfrak j\) of \(\mathfrak b\) acts trivially on
\(\mathfrak n\), then it is a direct Lie algebra summand of
\(\mathfrak g\).
\end{lemma}
\begin{proof}
The nilradical of \(\mathfrak r\) is characteristic and hence
is a nilpotent ideal of \(\mathfrak g\). Every nilpotent ideal of
\(\mathfrak g\) lies in \(\mathfrak r\), so the two nilradicals
agree. Triangularize the adjoint representation of
\(\mathfrak r_{\CC}\) on itself using Lie's theorem, and let
\(\chi_i\) be its diagonal characters. For \(X\in\mathfrak g\),
the automorphism \(e^{t\ad_X}|_{\mathfrak r}\) preserves the
finite multiset of characters by conjugacy of the adjoint
representations. Thus \(\chi_i\circ e^{t\ad_X}=\chi_i\) for
each \(i\), by continuity. Differentiating shows that
\(\ad_{[X,Y]}|_{\mathfrak r}\) is strictly upper triangular
for \(Y\in\mathfrak r\). Engel's theorem applied to the
ideal \([\mathfrak g,\mathfrak r]\) gives the inclusion.

If \([\mathfrak j,\mathfrak n]=0\), then
\([[\mathfrak j,\mathfrak j],\mathfrak r]=0\) by Jacobi and
the inclusion just proved. Since \(\mathfrak j\) is
semisimple, it commutes with \(\mathfrak r\) and with its
complementary ideal in \(\mathfrak b\).
\end{proof}

\begin{lemma}\label{lem:compact-obstruction}
A Lie algebra with a nonzero compact semisimple ideal admits
no negative Ricci inner product. For a compact Levi ideal
\(\mathfrak u\) of \(\mathfrak b\), its action on
\(\mathfrak n\) must therefore be faithful in the presence of
a negative Ricci inner product. If \(\mathfrak u\) is the sum
of all compact simple ideals of \(\mathfrak b\), faithfulness
is equivalent to the absence of a nonzero compact semisimple
ideal of \(\mathfrak g\).
\end{lemma}
\begin{proof}
The obstruction is \cite[Lemma~1.6]{JP} and
\cite[Proposition~2.1]{AraujoObstruction}. Indeed, on a connected
group with a left-invariant metric, a nonzero vector \(X\) in a
compact semisimple ideal generates a Killing field for left
translations. Its squared length at \(g\) is
\(\|\operatorname{Ad}(g^{-1})X\|^2\). The adjoint orbit of \(X\)
is compact: the connected adjoint action on the ideal is the
compact group of its inner automorphisms. The length attains its
maximum, contradicting the Killing field identity
\(\frac12\Delta |X^*|^2=|\nabla X^*|^2-\Ric(X^*,X^*)\)
under strictly negative Ricci curvature.

The kernel of the compact Levi action is a compact
semisimple ideal of \(\mathfrak b\). Lemma~\ref{lem:radical}
makes it an ideal of \(\mathfrak g\). Conversely, a compact
semisimple ideal \(\mathfrak j\) of \(\mathfrak g\) projects
injectively to a compact ideal of \(\mathfrak b\), since its
intersection with the solvable radical is zero. For an element
\(X+v\) of \(\mathfrak j\), with \(X\in\mathfrak u\), one has
\(v\in\mathfrak n\): projection of the radical component to
\(\mathfrak r/\mathfrak n\) is a homomorphism into an abelian
Lie algebra, and \(\mathfrak j\) is perfect. Also
\([\mathfrak j,\mathfrak r]=0\). Thus the action of \(X\) on
\(\mathfrak n\) equals \(-\ad_v|_{\mathfrak n}\). It is
nilpotent, and it is skew for a compact-invariant inner
product on \(\mathfrak n\), so it is zero. This proves the
last assertion.
\end{proof}

\begin{lemma}\label{lem:abelian-trace}
If \(W\) is an abelian ideal of a metric Lie algebra, then
\[
 \tr_W\Ric
 =\frac14\sum_{i,j}\|\operatorname{pr}_W[F_i,F_j]\|^2
  -\tr(\ad_H|_W),
\]
where the \(F_i\) form an orthonormal basis of \(W^\perp\) and
\(\tr_W\Ric\) is the trace of the compression to \(W\).
In particular, for a negative Ricci inner product and
\(W\ne0\), one has \(\tr(\ad_H|_W)>0\).
\end{lemma}
\begin{proof}
For \(w\in W\), the image of \(\ad_w\) lies in \(W\) and
\(\ad_w|_W=0\), so the Killing form vanishes on \(W\times W\).
Use the bases of \(W\) and \(W^\perp\) in
\eqref{eq:general-ricci}. The first term traced on \(W\) is
\(-\frac12\sum_{i,a}\|[F_i,w_a]\|^2\). It cancels the
contribution in the second term from pairs with one entry in
\(W\) and one in \(W^\perp\). The drift contributes
\(-\tr(\ad_H|_W)\).
\end{proof}

A central vector \(z\) satisfies
\(\Ric(z,z)=\frac14\sum_{i,j}\langle[e_i,e_j],z\rangle^2\geq0\).
The obstruction for solvable unimodular Lie groups goes back
to Dotti de Miatello \cite{DottiSolvable}. The semisimplicity
conclusion for unimodular groups was obtained by Dotti
Miatello, Leite and Miatello \cite{DLM}. B\"ohm and Lafuente
prove a version for cocompact isometric actions
\cite[Theorem~4.1]{BLUnimodular}. The abelian-ideal formula
gives a short proof for left-invariant metrics.

Also, if \(\mathfrak r\ne0\), the last nonzero term in the
derived series of \(\mathfrak r\) is a nonzero abelian ideal
of \(\mathfrak g\). Lemma~\ref{lem:abelian-trace} shows that
a negative Ricci metric in this case requires
nonunimodularity. In particular \(\mathfrak r\ne\mathfrak n\),
since a nilpotent radical has traceless adjoint action and
the trace character vanishes on the perfect Levi algebra.

We will also use the following full block formulas. On
\[
  (\mathfrak b\oplus\mathfrak a)\ltimes V
\]
with \(\mathfrak b\) semisimple and \(\mathfrak a,V\) abelian,
write the actions on \(V\) as \(\pi(X)\) and \(D_A\).
Use an inner product with the three summands orthogonal and
orthonormal bases \(X_i,A_j\) of the first two. The mean-curvature
vector lies in \(\mathfrak a\) and satisfies
\(\langle H,A\rangle=\tr D_A\). The Ricci blocks are
\begin{align}
 \Ric(X,X')&=\Ric_{\mathfrak b}(X,X')
               -\tr(\Sop(\pi(X))\Sop(\pi(X'))),\notag\\
 \Ric(A,A')&=-\tr(\Sop(D_A)\Sop(D_{A'})),\notag\\
 \Ric(X,A)&=-\tr(\Sop(\pi(X))\Sop(D_A)),\notag\\
 \Ric|_V&=\frac12\sum_i[\pi(X_i),\pi(X_i)^*]
         +\frac12\sum_j[D_{A_j},D_{A_j}^*]-\Sop(D_H).
 \label{eq:abelian-blocks}
\end{align}
The mixed terms with \(V\) vanish. For two first-block inputs,
the additional contributions come from the first and Killing
terms in \eqref{eq:general-ricci}. For two action operators
\(T,U\), they sum to
\(-\frac12\tr(T^*U+TU)=-\tr(\Sop(T)\Sop(U))\).
On \(V\), each action contributes
\(-\tfrac12 T^*T+\tfrac12 TT^*\); the Killing form vanishes
on \(V\times V\). The map that negates \(V\) and fixes
\(\mathfrak b\oplus\mathfrak a\) is an orthogonal automorphism
and gives the remaining mixed-block assertion.
\section{Attaching a split three-dimensional simple factor}
\label{sec:attachment}

\begin{theorem}
\label{thm:attachment}
Let \(\mathfrak h\) be a nonzero finite-dimensional real Lie algebra admitting a positive-definite inner product with strictly negative Ricci tensor. Then \(\mathfrak{sl}(2,\mathbb R)\oplus\mathfrak h\) admits such an inner product. The inner product can be chosen to restrict to the given one on \(\mathfrak h\). If \(\mathfrak h\) has nonzero nilradical, a single cross term between \(\mathfrak{sl}(2,\mathbb R)\) and that nilradical suffices.
\end{theorem}

\begin{proof}
Write \(R_{\mathfrak h}\) for the Ricci endomorphism of the given inner product, and let \(H_{\mathfrak h}\) satisfy \(\langle H_{\mathfrak h},w\rangle=\operatorname{tr}(\operatorname{ad}_{\mathfrak h}w)\). Since a one-dimensional Lie algebra is flat, \(\dim\mathfrak h\geq 2\). Choose a unit vector \(v\) in the kernel of \(w\mapsto\operatorname{tr}(\operatorname{ad}_{\mathfrak h}w)\). If the actual nilradical is nonzero, \(v\) can be chosen there. Indeed, for \(v\) in the nilradical, \(\operatorname{ad}_{\mathfrak h}v\) maps \(\mathfrak h\) into that ideal and acts nilpotently on it.

Fix a basis \(k_0,p_0,q_0\) of \(\mathfrak{sl}(2,\mathbb R)\) with
\[
 [k_0,p_0]=q_0,\qquad [k_0,q_0]=-p_0,\qquad [p_0,q_0]=-k_0.
\]
For \(t>0\), put \(a=t/2\) and set \(k=k_0\), \(p=\sqrt a\,p_0\), \(q=\sqrt a\,q_0\), and \(P=p+tv\). Declare \(k,P,q\) and an orthonormal basis \((e_i)\) of \(\mathfrak h\) containing \(v\) to be orthonormal for the new inner product. In the original direct sum, its only cross term in the basis \(k,p,q,(e_i)\) is \(\langle p,v\rangle=-t\). With \(L=\operatorname{ad}_{\mathfrak h}v\), the brackets in the orthonormal basis, apart from the original brackets in \(\mathfrak h\), are
\begin{equation}\label{eq:attachment-brackets}
 \begin{aligned}
 {}[k,P]&=q,& [k,q]&=-P+tv,\\
 [P,q]&=-ak,& [P,w]&=tLw \quad (w\in\mathfrak h).
 \end{aligned}
\end{equation}
The mean-curvature vector for the new inner product is \(H_{\mathfrak h}\). This follows from the trace-zero condition on \(v\), unimodularity of \(\mathfrak{sl}(2,\mathbb R)\), and the direct-sum bracket.

In the full Ricci formula \eqref{eq:general-ricci}, denote the negative and positive bracket-square summands by \(M_-\) and \(M_+\), respectively. Their entries from \eqref{eq:attachment-brackets}, together with those of the Killing form \(B\), are
\[
\begin{array}{c|ccc}
 & M_- & M_+ & B\\ \hline
 (k,k)&-1-t^2/2&a^2/2&-2\\
 (q,q)&-(1+a^2+t^2)/2&1/2&2a\\
 (P,P)&-(1+a^2+t^2\|L\|_{\mathrm{HS}}^2)/2&1/2&
          2a+t^2 B_{\mathfrak h}(v,v)\\
 (P,w)&-\frac t2\sum_i\langle Le_i,[w,e_i]\rangle&
       -\frac t2\langle v,w\rangle&tB_{\mathfrak h}(v,w)
\end{array}
\]
for \(w\in\mathfrak h\). The only bracket with both a \(P\)-component and an \(\mathfrak h\)-component is \([k,q]=-P+tv\), giving the indicated entry of \(M_+\). For \(w,z\in\mathfrak h\), the \(M_-\) contribution beyond the original one is \(-\tfrac12t^2\langle Lw,Lz\rangle\), from the index \(P\). The \(M_+\) contributions beyond the original one come from \([k,q]\) and \([P,e_i]\), and give
\(\tfrac12t^2(\langle v,w\rangle\langle v,z\rangle + \langle L^*w,L^*z\rangle)\).
The Killing form splits on the two original ideals. Its restriction to \(\mathfrak h\) is \(B_{\mathfrak h}\), and the displayed entries follow by evaluating it on \(P=p+tv\). In the mean-curvature term, \([H_{\mathfrak h},k]=[H_{\mathfrak h},q]=0\), \([H_{\mathfrak h},P]=t[H_{\mathfrak h},v]\), and the restriction to \(\mathfrak h\) is unchanged. Finally, for each fixed summation index in \(M_-\), a bracket involving \(k\) or \(q\) has output orthogonal to that for the other input in each mixed entry with \(k\) or \(q\). The bracket \([P,q]\) has output only in \(\mathbb Rk\) and \([k,P]\) only in \(\mathbb Rq\), so those mixed entries of \(M_+\) vanish as well. The Killing and mean-curvature terms also vanish in those entries.

Consequently \(k,q\) are Ricci eigenvectors with eigenvalues
\[
 \operatorname{Ric}(k,k)=-3t^2/8,\qquad
 \operatorname{Ric}(q,q)=-t/2-5t^2/8.
\]
On \(\mathbb RP\oplus\mathfrak h\) write the remaining block as
\(\begin{pmatrix}d_t&t\ell\\t\ell^*&R_t\end{pmatrix}\), where
\begin{align}
 d_t&=-t/2-t^2/8-t^2\|\Sop(L)\|_{\mathrm{HS}}^2,\notag\\
 \ell(w)&=-\tfrac12\left(
     \sum_i\langle Le_i,[w,e_i]\rangle+\langle v,w\rangle
     +B_{\mathfrak h}(v,w)+\langle[H_{\mathfrak h},v],w\rangle\right),\notag\\
 R_t&=R_{\mathfrak h}+\tfrac12t^2\bigl([L,L^*]+v\otimes v\bigr).
 \label{eq:attachment-remaining-block}
\end{align}
The expression for \(d_t\) uses
\(\|L\|_{\mathrm{HS}}^2+\tr(L^2)=2\|\Sop(L)\|_{\mathrm{HS}}^2\).
Here \(v\otimes v\) denotes \(w\mapsto\langle v,w\rangle v\).
Choose \(\lambda>0\) with \(R_{\mathfrak h}\leq-\lambda I\).
The last identity gives \(R_t\leq-\lambda I/2\) for all
sufficiently small \(t>0\). For such \(t\), the scalar Schur
complement satisfies
\[
 d_t-t^2\ell R_t^{-1}\ell^*
 \leq -t/2+2\|\ell\|^2t^2/\lambda<0
\]
after decreasing \(t\).
This proves negativity of the remaining block and hence of the full Ricci tensor.
\end{proof}

\begin{corollary}
\label{cor:attachment-iterate}
If \(\mathfrak h\ne0\) admits a strictly negative Ricci inner product, then \(\mathfrak{sl}(2,\mathbb R)^{\oplus m}\oplus\mathfrak h\) admits one for every integer \(m\geq 1\). In particular, if \(\mathfrak r\) is solvable with actual nonzero nilradical \(\mathfrak n\) and there is \(A\in\mathfrak r\) with
\(\operatorname{Re}\operatorname{Spec}(\operatorname{ad} A|_{\mathfrak n})>0\), then \(\mathfrak{sl}(2,\mathbb R)^{\oplus m}\oplus\mathfrak r\) admits one.
\end{corollary}

\begin{proof}
Apply Theorem~\ref{thm:attachment} successively. For the second assertion, the positive-stable hypothesis gives the initial metric on \(\mathfrak r\) by \cite[Theorem~2(2)]{NN}. The nilradical of each direct sum is the original \(\mathfrak n\), since its projection to the semisimple summand is a nilpotent ideal and therefore zero.
\end{proof}

\subsection{A semisimple counterexample to inheritance by the factors}

\begin{corollary}
\label{cor:attachment-semisimple}
For all integers \(m\geq1\) and \(d\geq3\), the semisimple Lie algebra
\[
  \mathfrak{sl}(2,\mathbb R)^{\oplus m}\oplus\mathfrak{sl}(d,\mathbb R)
\]
admits a strictly negative Ricci inner product. Each of its \(\mathfrak{sl}(2,\mathbb R)\) simple ideals admits no such inner product. In particular this applies to \(\mathfrak{sl}(2,\mathbb R)\oplus\mathfrak{sl}(3,\mathbb R)\).
\end{corollary}

\begin{proof}
Leite and Dotti de Miatello construct a strictly negative Ricci metric on \(\mathrm{SL}(d,\mathbb R)\) for \(d\geq3\) \cite{DL}. Apply Corollary~\ref{cor:attachment-iterate} to its Lie algebra. For an arbitrary inner product on \(\mathfrak{sl}(2,\mathbb R)\), the three-dimensional unimodular normal form gives cyclic brackets \([E_2,E_3]=-b_1E_1\), \([E_3,E_1]=b_2E_2\), and \([E_1,E_2]=b_3E_3\) in a suitable orthonormal basis, with each \(b_i>0\). The Ricci eigenvalues are
\[
 \tfrac12\bigl(b_1^2-(b_2-b_3)^2\bigr),\qquad
 \tfrac12\bigl(b_2^2-(b_1+b_3)^2\bigr),\qquad
 \tfrac12\bigl(b_3^2-(b_1+b_2)^2\bigr)
\]
by \cite[Theorem~4.3 and Corollary~4.7]{Milnor}. The first is nonnegative when \(|b_2-b_3|\le b_1\). In the other two cases \(b_2>b_1+b_3\) or \(b_3>b_1+b_2\), the second or third eigenvalue is positive, respectively.
\end{proof}

\subsection{Reverse inheritance for particular solvable factors}

In the following cases, the abelian ideal trace identity also determines when the solvable factor of the direct sum has negative Ricci curvature.

\begin{theorem}
\label{thm:attachment-reverse}
Let \(m\geq1\) and let \(\mathfrak r\) be a real solvable Lie algebra with actual nonzero nilradical \(\mathfrak n\). Suppose \(\mathfrak n\) is abelian, Heisenberg, or standard filiform. Then
\(\mathfrak{sl}(2,\mathbb R)^{\oplus m}\oplus\mathfrak r\) admits a strictly negative Ricci inner product if and only if \(\mathfrak r\) admits one. More explicitly, the respective equivalent conditions are as follows.
\begin{enumerate}
 \item If \(\mathfrak n\) is abelian, there is \(A\in\mathfrak r\) with
 \(\operatorname{Re}\operatorname{Spec}(\operatorname{ad} A|_{\mathfrak n})>0\).
 \item If \(\mathfrak n\) is Heisenberg, let \(\lambda(A)\) be the scalar of \(\operatorname{ad}A\) on the center \(\mathfrak z\), and let \(d_j(A)\) be the eigenvalues of the induced operator on \(\mathfrak n/\mathfrak z\), with algebraic multiplicity. There is \(A\in\mathfrak r\) with
 \[
 \lambda(A)+\sum_{\operatorname{Re}d_j(A)<0}\operatorname{Re}d_j(A)>0.
 \]
 \item If \(\mathfrak n\) is the standard filiform algebra of dimension \(l\geq4\), use a basis \(X_1,\ldots,X_l\) with \([X_1,X_j]=X_{j+1}\) for \(2\leq j<l\) and no other independent nonzero bracket. Write \(\mathfrak z=\mathbb RX_l\), \(I=\operatorname{span}\{X_2,\ldots,X_l\}\), and let \(\lambda(A)\) be the scalar of \(\operatorname{ad}A\) on \(\mathfrak z\). There is \(A\in\mathfrak r\) with
 \[
  \lambda(A)>0,\qquad \operatorname{tr}(\operatorname{ad}A|_I)>0.
 \]
\end{enumerate}
\end{theorem}

\begin{proof}
Each listed condition implies that \(\mathfrak r\) admits a strictly negative Ricci inner product by the abelian corollary to Theorem~2, Theorem~3, and Theorem~4 of \cite{NN}, respectively. The resulting metric extends by Corollary~\ref{cor:attachment-iterate}.

For necessity, give \(\mathfrak g=\mathfrak{sl}(2,\mathbb R)^{\oplus m}\oplus\mathfrak r\) any strictly negative Ricci inner product. Its actual nilradical is \(\mathfrak n\). Let \(H\) be the mean-curvature vector and set \(A=\operatorname{pr}_{\mathfrak r}H\) using the algebraic direct sum. Write \(D_X=\operatorname{ad}X|_{\mathfrak n}\) for \(X\in\mathfrak g\). The simple summands centralize \(\mathfrak n\), so \(D_H=D_A\). By Lemma~\ref{lem:abelian-trace}, every nonzero abelian ideal \(W\) of \(\mathfrak g\) contained in \(\mathfrak n\) satisfies
\begin{equation}
\label{eq:attachment-trace-test}
 \operatorname{tr}(D_H|_W)>0.
\end{equation}
The induced \(D_X\) commute on the center \(\mathfrak z=\mathfrak z(\mathfrak n)\). Indeed, the \(D_X\) depend only on the radical projections, \([\mathfrak r,\mathfrak r]\subseteq\mathfrak n\), and each inner derivation of \(\mathfrak n\) acts trivially on \(\mathfrak z\). If \(D_H|_{\mathfrak z}\) had an eigenvalue with nonpositive real part, the sum of its complex generalized eigenspaces with nonpositive real parts would define a nonzero real subspace \(W\) of \(\mathfrak z\) invariant under every \(D_X\). This gives an abelian ideal of \(\mathfrak g\) whose trace for \(D_H\) is nonpositive, contrary to \eqref{eq:attachment-trace-test}. Therefore \(D_H|_{\mathfrak z}\) is positive stable. This proves the abelian case.

For a Heisenberg nilradical, write \([x,y]=\omega(\bar x,\bar y)Z\) with \(0\ne Z\in\mathfrak z\), and let \(\overline D_H\) be induced by \(D_H\) on \(\mathfrak n/\mathfrak z\). The induced \(D_X\) commute on that quotient as well, since inner derivations of \(\mathfrak n\) have image in \(\mathfrak z\). Let \(V_-\) be the real subspace of the quotient obtained from generalized eigenspaces of \(\overline D_H\) with negative real parts. Its inverse image \(W\) in \(\mathfrak n\) contains \(\mathfrak z\) and is invariant under every \(D_X\), so it is an ideal of \(\mathfrak g\). It is abelian. In fact, for \(x,y\in V_-\), the derivation identity gives
\[
 \omega(e^{t\overline D_H}x,e^{t\overline D_H}y)
 =e^{t\lambda(A)}\omega(x,y).
\]
The left hand side tends to zero as \(t\to+\infty\), whereas \(\lambda(A)>0\) by the center argument. Hence \(\omega(x,y)=0\). Use \eqref{eq:attachment-trace-test} on \(W\) to obtain
\(\lambda(A)+\sum_{\operatorname{Re}d_j(A)<0}\operatorname{Re}d_j(A)>0\).

For a standard filiform nilradical, \(I\) is a characteristic abelian ideal. To see uniqueness among codimension-one abelian ideals, any codimension-one ideal contains \([\mathfrak n,\mathfrak n]\), in particular \(X_3\). An element with nonzero \(X_1\)-component has nonzero bracket with \(X_3\), so an abelian ideal of that dimension must equal \(I\). Every derivation preserves \(I\), since its exponential is an automorphism. Thus \eqref{eq:attachment-trace-test} yields the trace inequality on \(I\), while the center argument supplies \(\lambda(A)>0\). This completes the three cases.
\end{proof}
\section{Proof of the stable radical criterion}\label{sec:mixed-proof}

Necessity in Theorem~\ref{thm:mixed} follows from
Lemma~\ref{lem:compact-obstruction}. For sufficiency, first remove
the split simple ideals that act trivially on \(\mathfrak n\).
They are direct summands by Lemma~\ref{lem:radical} and will be
adjoined by Theorem~\ref{thm:attachment}. The remaining Levi
algebra is denoted by \(\mathfrak b=\mathfrak u\oplus
\mathfrak l_-\oplus\mathfrak s\), where \(\mathfrak s\) is the
sum of the split factors that act nontrivially. The radical,
the actual nilradical, and the stable direction are the same
after this removal.

\subsection{Reduction to joint characters}

Two Lie brackets on the same vector space are isomorphic if
one is obtained from the other by a change of basis. If a
bracket is a limit of brackets isomorphic to a given bracket,
and the limit admits a negative Ricci inner product, so does
the given bracket. Indeed, for the fixed inner product on
the vector space all entries of \eqref{eq:general-ricci}
depend continuously on the bracket. Negative definiteness is
open, and an isomorphism pulls the inner product back
\cite[Proposition~1]{NN}.

\begin{lemma}\label{lem:stable-reduction}
There is an abelian complementary subalgebra
\(\mathfrak a\) to \(\mathfrak n\) in \(\mathfrak r\) that
centralizes \(\mathfrak b\) and contains a stable direction
\(Y_0\). The bracket of \(\mathfrak g\) degenerates in two
steps to
\begin{equation}\label{eq:reduced}
 (\mathfrak b\oplus\mathfrak a)\ltimes V,
 \qquad V=\bigoplus_E E\ \text{abelian}.
\end{equation}
On \(V\) the Levi representation \(\pi\) is the original
representation on the vector space \(\mathfrak n\) and
preserves each \(E\). The \(\mathfrak a\)-action has the form
\begin{equation}\label{eq:characters}
 D_A|_E=\alpha_E(A)I+\beta_E(A)J_E ,
\end{equation}
with real linear characters \(\alpha_E,\beta_E\). On a real
character block the second term is zero. On a nonreal block
\(J_E^2=-I\) and \(J_E\) commutes with \(\pi(\mathfrak b)\).
The joint map \(A\mapsto(\alpha_E(A),\beta_E(A))_E\) is
injective and \(\alpha_E(Y_0)>0\) for every \(E\).
\end{lemma}
\begin{proof}
The quotient \(\mathfrak r/\mathfrak n\) is trivial as a
\(\mathfrak b\)-module by Lemma~\ref{lem:radical}. Complete
reducibility gives a representative \(Y_0\) of the class
of \(Y\) that centralizes \(\mathfrak b\). On the lower
central series of \(\mathfrak n\), \(\ad_v|_{\mathfrak n}\)
raises filtration degree for \(v\in\mathfrak n\), while any
derivation preserves the filtration. Thus
\(\ad_{Y_0}|_{\mathfrak n}\) and \(\ad_Y|_{\mathfrak n}\)
have the same characteristic polynomial. Write
\(D_{Y_0}=\ad_{Y_0}|_{\mathfrak n}\). It is invertible.
The centralizer \(\mathfrak a=C_{\mathfrak r}(Y_0)\) is
a complement: the image of \(A\in\mathfrak r\) in
\(\mathfrak r/\mathfrak n\) has a unique representative
\(A-D_{Y_0}^{-1}[Y_0,A]\) in the centralizer. Jacobi gives
\([\mathfrak a,\mathfrak a]\subseteq\mathfrak n\cap
C_{\mathfrak r}(Y_0)=0\) and
\([\mathfrak b,\mathfrak a]\subseteq\mathfrak n\cap
C_{\mathfrak r}(Y_0)=0\).

For a linear automorphism \(h\) write the transported bracket
as \(h\mu(h^{-1}\,\cdot,h^{-1}\,\cdot)\). Scaling
\(\mathfrak n\) by a factor tending to infinity and fixing
\(\mathfrak b\oplus\mathfrak a\) makes \([\mathfrak n,\mathfrak n]\)
vanish and preserves all actions. On the resulting abelian
space, split the commuting family \(D_A\) into joint generalized
character spaces, pairing complex-conjugate characters over
\(\RR\). Every block is invariant under the commuting Levi action.
On a paired block \(J_E\) is multiplication by \(i\) on the
chosen character space and by \(-i\) on its conjugate.
Before the second degeneration, \eqref{eq:characters}
has an additional commuting nilpotent term \(N_A\) on each block.
For a basis of \(\mathfrak a\), the common flag of kernels of
all products of \(j\) of these nilpotent operators exhausts
the block, by simultaneous strict triangularization. Each
\(N_A\) maps level \(j\) into level \(j-1\). Split the flag
\(\mathfrak b\)-equivariantly by complete reducibility, using
complex-linear splittings on nonreal blocks. Scaling the
level-\(j\) summands by \(s^j\) as \(s\to\infty\) kills each
\(N_A\) and preserves the Levi and semisimple character
actions. Stability gives \(\alpha_E(Y_0)>0\).

If \(A\) is in the kernel of the joint character map,
its adjoint action on the original \(\mathfrak n\) is
nilpotent. The space \(\mathfrak n+\RR A\) is an ideal of the
original \(\mathfrak g\). Lie's theorem for the action of
\(\mathfrak r_\CC\) on \(\mathfrak n_\CC\) makes
\(\ad_{v+tA}|_{\mathfrak n}\) strictly upper triangular
for \(v\in\mathfrak n\) and \(t\in\RR\), since
\(\ad_v|_{\mathfrak n}\) and \(\ad_A|_{\mathfrak n}\)
both have zero diagonals in this triangularization. The
adjoint on the quotient of this ideal by \(\mathfrak n\)
is zero. Engel's theorem and the definition of the actual
nilradical imply \(A\in\mathfrak n\), hence \(A=0\).
\end{proof}
\subsection{The compact Levi action}\label{sec:compact-action}

\begin{lemma}\label{lem:boost}
Let $\pi:\uu\to\mathfrak{so}(E)$ be a compact semisimple representation
and let $\ttc\subset\uu$ be a Cartan subalgebra. If specified, let $J$
be an orthogonal complex structure commuting with $\pi(\uu)$.
There are positive self-adjoint $h_t$ commuting with $\pi(\ttc)$ and
with $J$ such that, for $R$ in the sum of the real root planes,
\[
 \|\Sop(h_t\pi(R)h_t^{-1})\|_{HS}^2
 \ge \sinh^2(t)\|\pi(R)\|_{HS}^2.
\]
For $H\in\ttc$, the operator $h_t\pi(H)h_t^{-1}=\pi(H)$ remains skew.
\end{lemma}
\begin{proof}
Without a specified $J$, decompose the commuting skew action of $\ttc$
orthogonally into one-dimensional zero-weight lines and two-dimensional
rotation planes. With $J$ specified, choose a unitary common weight basis
of the complex space $(E,J)$ and take its complex lines as real
two-dimensional blocks. Write this decomposition as $E=\bigoplus_j E_j$.
Choose distinct integers $c_j$ and put $h_t|_{E_j}=e^{c_jt}I$.

The diagonal compression $X\mapsto\pi(X)_{jj}$ is $\ttc$-equivariant
and takes values in $\mathfrak{so}(E_j)$, because the original action is
skew. This target is zero in dimension one and is a trivial $\ttc$-module
in dimension two. Every real root plane is a nontrivial $\ttc$-module.
Consequently $\pi(R)_{jj}=0$.
For $i\ne j$, skewness gives
\[
 \Sop(h_t\pi(R)h_t^{-1})_{ij}
       =\sinh((c_i-c_j)t)\pi(R)_{ij}.
\]
The diagonal blocks vanish. Summing their squared Hilbert--Schmidt norms
proves the estimate. The construction commutes with the Cartan action,
and in the complex case each $E_j$ is a $J$-line.
\end{proof}

\begin{lemma}\label{lem:cross}
Suppose $B_Z$ commutes with the fixed skew Cartan action and
$A_X=h_t\pi(X)h_t^{-1}$. Then
\[
 \tr\bigl(\Sop(A_X)\Sop(B_Z)\bigr)=0\qquad(X\in\uu).
\]
In particular this applies to a commuting representation of any other
semisimple factor transported through the same $h_t$.
\end{lemma}
\begin{proof}
The transpose of $B_Z$ also commutes with the skew Cartan, so the same
holds for $\Sop(B_Z)$. Taking symmetric parts intertwines commutators
with a skew operator. Thus $X\mapsto\Sop(A_X)$ is Cartan-equivariant,
vanishes on $\ttc$, and sends root planes to nonzero-weight modules.
These modules are orthogonal to the zero-weight space containing
$\Sop(B_Z)$. Equivalently, use
\[
 \tr([\pi(H),\Sop(A_{Y_\alpha})]\Sop(B_Z))
 =\tr(\Sop(A_{Y_\alpha})[\Sop(B_Z),\pi(H)])=0
\]
with $\alpha(H)\ne0$ and the root relation.
\end{proof}

To choose the intrinsic metric on \(\uu\), put \(Q=-B_\uu\).
For each positive root choose a \(Q\)-orthonormal pair \(X_\alpha,Y_\alpha\)
with \([H,X_\alpha]=\alpha(H)Y_\alpha\) and
\([H,Y_\alpha]=-\alpha(H)X_\alpha\) for \(H\in\ttc\).
Keep \(Q\) on \(\ttc\). Give \(X_\alpha,Y_\alpha\) lengths one and
two, respectively, and keep these directions orthogonal. In
\eqref{eq:general-ricci} the first sum at \((H,H)\) is
\(-\frac{17}{8}\sum_{\alpha>0}\alpha(H)^2\), since
\(X_\alpha,Y_\alpha/2\) are unit vectors. The Killing term is
\(\sum_{\alpha>0}\alpha(H)^2\) by the root decomposition.
The positive bracket sum is \(\frac18\sum_{\alpha>0}\alpha(H)^2\):
the only contributing ordered pairs are \(X_\alpha,Y_\alpha/2\)
and their reversals, and \(Q([X_\alpha,Y_\alpha],H)=\alpha(H)\).
The mean-curvature term is zero by semisimplicity. Therefore
\[
 \Ric_{\uu}(H,H)=-\sum_{\alpha>0}\alpha(H)^2=-Q(H,H)/2 .
\]
Its fixed Cartan/root cross terms satisfy a Young estimate
\begin{equation}\label{eq:intrinsic-bound}
 \Ric_{\uu}(H+R,H+R)
 \le -c_0\|H\|_Q^2+C_0\|R\|_Q^2,
 \qquad c_0>0.
\end{equation}
On each character block take a compact-invariant reference inner product
with $J_E$ orthogonal when present. Lemma~\ref{lem:boost} then gives a
deformation commuting with every $D_A$. For the actual choice of reference
metric one can first split off the real split representations as below,
and use this construction on their multiplicity spaces. Faithfulness
makes the total representation trace form positive definite. Hence, for
large fixed $t$,
\begin{equation}\label{eq:compact-negative}
 \Ric_{\uu}(X,X)-\|\Sop(h_t\pi(X)h_t^{-1})\|_{HS}^2<0
 \qquad(X\ne0).
\end{equation}
Lemma~\ref{lem:cross} makes the compact/noncompact Gram cross
exactly zero at this stage.

\subsection{The nontrivially acting split factors}\label{sec:split-action}

Every irreducible finite-dimensional real representation of
$\ssplit=\bigoplus_i\mathfrak{sl}(2,\RR)$ is
\[
 W_\lambda=\bigotimes_i\Sym^{m_{\lambda i}}\RR^2,
 \qquad m_{\lambda i}\ge0,
\]
and is absolutely real. One can see this by choosing real highest-weight
vectors and applying the real lowering operators; its complexification
is the irreducible module of the same highest weight. Consequently
$\End_{\ssplit}(W_\lambda)=\RR$ and
\begin{equation}\label{eq:isotypic}
 V=\bigoplus_\lambda W_\lambda\otimes_{\RR}M_\lambda.
\end{equation}
Every operator centralizing $\ssplit$ acts as $I\otimes B_\lambda$.
In particular, $\uu,\lln,\aaa$ act on the multiplicity spaces, and a
nonreal character structure is $J_E=I\otimes J_M$. The decomposition can be taken separately on each joint character block.

Give $\Sym^m\RR^2$ the metric inherited from the standard Euclidean
tensor power, and use product metrics on $W_\lambda\otimes M_\lambda$.
For
\[
 K=\begin{pmatrix}0&-1\\1&0\end{pmatrix},\quad
 P=\begin{pmatrix}0&1\\1&0\end{pmatrix},\quad
 H=\begin{pmatrix}1&0\\0&-1\end{pmatrix},
\]
their representation operators satisfy $\pi(K)^t=-\pi(K)$,
$\pi(P)^t=\pi(P)$ and $\pi(H)^t=\pi(H)$.
On the multiplicity spaces choose the compact-invariant reference
metrics and the deformation of Subsection~\ref{sec:compact-action}. This commutes with the split
action. Next apply
\[
 h_s=\exp\left(\sum_{i=1}^k s_i\pi(H_i)\right).
\]
It acts on the $W$ factors, commutes with the preceding deformation, and
does not change any multiplicity-space action.

The brackets are $[K,P]=-2H$, $[P,H]=2K$, $[H,K]=-2P$. Put
$x_i=\sinh(2s_i)$, $c_i=\cosh(2s_i)$ and
$q_i=\tr_V(\pi(P_i)^2)$. Nontriviality of the $i$th simple action gives
$q_i>0$. The compact $K_i$ rotation conjugates $H_i$ to $P_i$, so
\[
 \tr(\pi(H_i)^2)=q_i,\qquad \tr(\pi(P_i)\pi(H_i))=0.
\]
Direct conjugation gives
\[
 \Sop(A_{K_i})=-x_i\pi(P_i),\quad
 \Sop(A_{P_i})=c_i\pi(P_i),\quad
 \Sop(A_{H_i})=\pi(H_i).
\]
The intrinsic Ricci form for orthonormal $K_i,P_i,H_i$ is
$\operatorname{diag}(2,-6,-6)$. Subtracting the symmetric-action Gram
therefore gives the $K_i,P_i$ block
\begin{equation}\label{eq:split-block}
 \begin{pmatrix}
 2-q_ix_i^2&q_ix_ic_i\\
 q_ix_ic_i&-6-q_ic_i^2
 \end{pmatrix},
 \qquad \det=-12-2q_i+4q_ix_i^2,
\end{equation}
and the $H_i$ value is $-6-q_i$. Choose $s_i$ so that
$x_i^2>1/2+3/q_i$. The determinant is positive and the second diagonal
entry is negative, hence the whole block is strictly negative.

All remaining split-factor crosses vanish. Against a multiplicity-space
operator, the product trace factors as
\[
 \tr_W(\Sop(A_X))\,\tr_M(\Sop(B))=0.
\]
Against a different split factor, two distinct tensor factors each have
zero trace. These identities hold on each isotypic summand, including
summands where one action is trivial. The $\uu/\lln$ part is already
negative by \eqref{eq:compact-negative}, Lemma~\ref{lem:cross}, and the
given negative metric on $\lln$. The entire Levi block is now negative.

\subsection{Imaginary characters and the nilradical block}

The final perturbation detects pure imaginary character directions while
keeping the radical/base Ricci cross zero. On each nonreal multiplicity
block choose a symmetric orthogonal involution $T$ with $TJ=-JT$, and
put $h_\varepsilon=e^{\varepsilon T}$, $\varepsilon\ne0$ small.
On $W\otimes M$ use $I\otimes h_\varepsilon$; on real-character blocks
use the identity. With $L=TJ$, direct multiplication gives
\[
 h_\varepsilon^{-1}Lh_\varepsilon
   =\cosh(2\varepsilon)L-\sinh(2\varepsilon)J,
 \qquad
 \Sop(h_\varepsilon Jh_\varepsilon^{-1})
   =\sinh(2\varepsilon)L.
\]
For any complex-linear $B=\pi(X)$ on such a block,
$\tr(BL)=0$. Since the entire base $\bb$ is perfect,
$\tr B=\tr(BJ)=0$. Thus
\begin{equation}\label{eq:anti-cross}
 \tr\bigl(\Sop(h_\varepsilon Bh_\varepsilon^{-1})L\bigr)=0.
\end{equation}
Equation
\eqref{eq:anti-cross} and the scalar trace identity give exact zero for
the $\aaa/\bb$ symmetric Gram. Meanwhile
\[
 G(A,B)=\sum_E\dim_{\RR}(E)
 \left(\alpha_E(A)\alpha_E(B)
 +\sinh^2(2\varepsilon)\beta_E(A)\beta_E(B)\right)
\]
is positive definite by injectivity of the joint character map.

Fix the earlier parameters. Openness retains strict negativity of
the full Levi block for sufficiently small nonzero $\varepsilon$,
including its possibly nonzero compact/noncompact mixed terms.

Let $\tau=\tr D_{Y_0}>0$. At $\varepsilon=0$, all earlier deformations
commute with $D_{Y_0}$, which is normal with positive definite symmetric part. Consequently
\[
 C_{Y_0}=\tfrac12[D_{Y_0}',D_{Y_0}'{}^t]
                  -\tau\Sop(D_{Y_0}')<0
\]
also holds for sufficiently small nonzero $\varepsilon$.
Give $\aaa$ an inner product with $Y_0\perp\ker\theta$,
$|Y_0|=R^{-1}$, where $\theta(A)=\tr D_A$, and keep the metric on
$\ker\theta$ fixed. The mean-curvature vector is
$R^2\theta(Y_0)Y_0$.

With the parameters fixed in this order, \eqref{eq:abelian-blocks}
gives
\[
 \Ric|_V=R^2C_{Y_0}
 +\frac12\sum_j[\pi'(X_j),\pi'(X_j)^*]
 +\frac12\sum_m[D'_{B_m},D'_{B_m}{}^*],
\]
where \(X_j\) is a fixed orthonormal Levi basis and \(B_m\) a fixed
orthonormal basis of \(\ker\theta\). The sums are independent of \(R\),
so the \(V\) block is negative for large \(R\). The
\(\mathfrak a\) block is \(-G\), its cross with the Levi block
vanishes by \eqref{eq:anti-cross}, and all mixed terms with \(V\)
vanish by \eqref{eq:abelian-blocks}. Thus the limit bracket carries
a negative Ricci metric. Apply the bracket-limit argument at each
of the two degenerations in Lemma~\ref{lem:stable-reduction}, then
adjoin the split factors removed at the start. This completes
the proof of Theorem~\ref{thm:mixed}.

\begin{corollary}\label{thm:compact}
Suppose that the Levi algebra of \(\mathfrak g\) is compact and
semisimple, \(\mathfrak n\ne0\), and \(\mathfrak r\) contains a
stable direction. The following conditions are equivalent:
\(\mathfrak g\) admits a negative Ricci inner product, the Levi
action on \(\mathfrak n\) is faithful, and \(\mathfrak g\) has
no nonzero compact semisimple ideal.
\end{corollary}
\begin{proof}
Apply Theorem~\ref{thm:mixed} and Lemma~\ref{lem:compact-obstruction}.
\end{proof}

\begin{corollary}\label{thm:negative-base}
Suppose the real semisimple Levi algebra of \(\mathfrak g\)
has a given negative Ricci inner product, \(\mathfrak n\ne0\),
and \(\mathfrak r\) contains a stable direction. Then
\(\mathfrak g\) has a negative Ricci inner product restricting
to the given one on the Levi algebra.
\end{corollary}
\begin{proof}
Use the preceding proof with \(\mathfrak u=\mathfrak s=0\).
Both bracket changes fix the Levi algebra, and its given
inner product serves as the intrinsic metric throughout
\eqref{eq:abelian-blocks}.
\end{proof}

In comparison with \cite[Theorem~5.2]{WillII} (arXiv version~2),
Corollary~\ref{thm:negative-base} keeps the stable direction and
the negative Ricci metric on the Levi algebra. The proof here
permits a nonnormal original radical action, nonzero
directions with purely imaginary joint characters, and an
arbitrary negative Ricci metric on the Levi algebra. In the
compact case, \cite[Theorem~5.3]{LW} already permits a positive
real-diagonalizable derivation in a one-dimensional complement
subject to further conditions on the representation. The
argument above uses faithfulness of the total compact action
and permits a stable derivation with nonreal spectrum or
nontrivial Jordan blocks and a complement of higher dimension.
\section{Stable radical models and quantitative metrics}\label{sec:models}

\subsection{One-dimensional complements and Heintze radicals}

Let \(\mathfrak n\ne0\) be nilpotent and \(D\in\operatorname{Der}(\mathfrak n)\)
have all eigenvalues in the open right half-plane. If
\(\pi\colon\mathfrak b\to\operatorname{Der}(\mathfrak n)\)
commutes with \(D\), the Lie algebra
\[
 (\mathfrak b\oplus\RR Z)\ltimes\mathfrak n,\qquad
 [Z,v]=Dv,\quad [X,v]=\pi(X)v \quad(X\in\mathfrak b)
\]
has actual nilradical \(\mathfrak n\). Indeed its solvable radical
is \(\RR Z\ltimes\mathfrak n\); on \(\mathfrak n\), the adjoint of
\(tZ+v\) has the same characteristic polynomial as \(tD\) by the
lower central filtration. A nilpotent ideal cannot contain an
element with \(t\ne0\), since it would then contain
\([tZ+v,\mathfrak n]=\mathfrak n\) and have a nonnilpotent
adjoint action on itself. Theorem~\ref{thm:mixed} applies whenever
\(\mathfrak b\) has its stated form. Conversely, for a stable
radical with \(\dim\mathfrak r/\mathfrak n=1\),
Lemma~\ref{lem:stable-reduction} gives precisely this presentation.

Suppose \(\mathfrak r\) is the Lie algebra of a simply connected
solvable Lie group of dimension at least two carrying a
left-invariant metric of strictly negative sectional curvature.
Heintze's theorem gives
\(\mathfrak r=\RR Z\ltimes\mathfrak n\) with
\(D=\ad_Z|_{\mathfrak n}\) stable after orienting \(Z\)
\cite[Theorem~3]{Heintze}. Here the derived algebra and the
actual nilradical coincide, since \(D\) is invertible.
In particular, for a compact Levi algebra acting on such a
radical, existence of a negative Ricci inner product is
equivalent to faithfulness on \(\mathfrak n\).

\subsection{A higher-dimensional complement with imaginary characters}

Let \(\mathfrak u=\mathfrak{su}(2)\) act in the standard way on
each copy of \(\CC^2\) in
\(V=(\CC^2\oplus\CC^2)_{\RR}\), and take \(V\) abelian. Let
\(\mathfrak a=\RR Z\oplus\RR T\) be abelian and commute with
\(\mathfrak u\). Write \(J\) for multiplication by \(i\) on
each copy and set
\[
 D_Z=(I+J)\oplus(I+2J),\qquad D_T=J\oplus J .
\]
For \(aZ+bT\) the real parts on both blocks equal \(a\), and
the imaginary coefficients are \(a+b\) and \(2a+b\).
Every nonzero vector in \(\mathfrak a\) has at least one
nonzero joint character. It follows that \(V\) is the actual
nilradical: for any nonzero \(\mathfrak a\)-component an
adjoint on \(V\) has a nonzero eigenvalue, so a nilpotent
ideal containing that element would contain the image of
that adjoint on \(V\) and would have a nonnilpotent
adjoint on its intersection with \(V\). The vector \(Z\)
is stable, \(T\) has purely imaginary spectrum, and a
stable vector cannot have both imaginary coefficients zero.
Theorem~\ref{thm:mixed} yields a negative Ricci inner
product on this thirteen-dimensional algebra.

\subsection{An explicit scalar extension of any faithful compact module}

\begin{proposition}\label{prop:scalar-construction}
Let \(\mathfrak u\ne0\) be compact semisimple and let
\(\pi\colon\mathfrak u\to\End(V)\) be any faithful real
representation, with \(d=\dim_{\RR}V\). Fix \(c>0\), make
\(V\) abelian, and let \(Z\) act on \(V\) by \(cI\) and
commute with \(\mathfrak u\). Choose an invariant inner
product \(q_0\) on \(V\) and any \(q_0\)-orthonormal basis.
Put
\[
 Q(X,Y)=-\tr(\pi(X)\pi(Y)),\qquad
 h=\operatorname{diag}(1,d,\ldots,d^{d-1}).
\]
On \((\mathfrak u\oplus\RR Z)\ltimes V\), make the three
summands orthogonal and set
\[
 g|_{\mathfrak u}=Q,\qquad
 g(v,w)=q_0(hv,hw),\qquad |Z|=c/d^d .
\]
Then
\[
 \Ric_g\leq-\frac{d^2-d+d^{-2}}4\,g .
\]
\end{proposition}
\begin{proof}
Faithfulness gives \(d\geq3\). The representation identifies
\((\mathfrak u,Q)\) with a subalgebra of \(\mathfrak{so}(d)\)
for the Hilbert--Schmidt inner product. The action
\(\ad_{\pi(X)}\) on \(\mathfrak{so}(d)\) is skew and
preserves the subalgebra and its orthogonal complement, so
\[
 -B_{\mathfrak u}(X,X)
 \leq -B_{\mathfrak{so}(d)}(\pi(X),\pi(X))
 =(d-2)Q(X,X).
\]
The last equality follows by identifying
\(\mathfrak{so}(d)\) with \(\Lambda^2V\): the trace of
\((T\otimes I+I\otimes T)^2\) on \(\Lambda^2V\) equals
\((d-2)\tr(T^2)+(\tr T)^2\), by subtracting the
flip-weighted trace on \(V\otimes V\) and dividing by two.
Here \(T=\pi(X)\) is skew. Since \(Q\) is bi-invariant,
\(\Ric_{\mathfrak u}=-B_{\mathfrak u}/4\).

In the metric basis \(h^{-1}e_i\) of \(V\), the representation
is \(A_X=h\pi(X)h^{-1}\). Its symmetric part has off-diagonal
entries
\[
 (\Sop(A_X))_{ij}
   =\frac{d^{i-j}-d^{j-i}}2\,\pi(X)_{ij}
\]
and zero diagonal. Hence
\(\|\Sop(A_X)\|_{\mathrm{HS}}^2
\geq (d-d^{-1})^2 Q(X,X)/4\). The compact Ricci block
in \eqref{eq:abelian-blocks} is therefore bounded above by
\(-\gamma_d Q\), for \(\gamma_d=(d^2-d+d^{-2})/4\).

For the unit vector \(e_0=(d^d/c)Z\) we have
\(\Ric(e_0,e_0)=-d^{2d+1}\) and \(\Ric(e_0,X)=0\)
by the scalar action and \(\tr A_X=0\). On \(V\), for a
\(Q\)-orthonormal basis \(X_1,\ldots,X_k\), formula
\eqref{eq:abelian-blocks} gives
\[
 \Ric|_V=\tfrac12\sum_{j=1}^k[A_{X_j},A_{X_j}^*]-d^{2d+1}I.
\]
Each commutator has half its operator norm at most
\(\|A_{X_j}\|_{\mathrm{HS}}^2\leq d^{2d-2}\).
Since \(k=\dim\mathfrak u\leq d(d-1)/2\), this yields
\(\Ric|_V\leq-d^{2d-2}(d^3-k)I\).
Both \(d^{2d+1}\) and \(d^{2d-2}(d^3-k)\) exceed \(\gamma_d\)
for \(d\geq3\). All mixed terms with \(V\) vanish by
\eqref{eq:abelian-blocks}.
\end{proof}

\begin{remark}
Proposition~\ref{prop:scalar-construction} applies to the small
complex modules discussed in \cite[Remark~4.3]{LW}, regarded
as real modules. For the fundamental \(A_2\) module, the vector
and spin \(B_2\) modules, and the fundamental \(G_2\) module,
the total dimensions are \(15,21,19,29\), respectively.
In particular, the complex vector modules of \(B_2\) and \(G_2\)
give two copies of their real vector modules. The existence
conclusions for those two realifications can also be
obtained from \cite[Theorem~3.1]{LW}. For \(B_2\), take the
real Cartan plane of weight \(e_1\) in the first copy and
that of weight \(e_2\) in the second, using the vector weights
\(0,\pm e_1,\pm e_2\). For \(G_2\), take the planes of
weights \(\mu_1,\mu_2\) in different copies, using vector
weights \(0,\pm\mu_1,\pm\mu_2,\pm\mu_3\) with
\(\mu_1+\mu_2+\mu_3=0\). In either case, no doubled nonzero
weight is a root, so each root plane maps the sum of these
selected planes into its orthogonal complement. The \(B_2\)
long roots touch both weight pairs and each short root
touches its own pair. The \(G_2\) long root
\(\mu_i-\mu_j\) touches the \(i\)th and \(j\)th pairs,
and each short root touches all three by the weight strings.
This gives a nonzero restriction for each root plane.
The Hilbert--Schmidt pairing on restrictions to the
selected planes is Cartan-invariant. Distinct real root
planes are orthogonal for this pairing, which is scalar
on each root plane, giving the restricted trace conditions
of that theorem. The single-copy models of dimensions
\(16\) and \(22\) in
Corollary~\ref{cor:small-compact-dimensions} use different real
modules.
\end{remark}

\begin{corollary}\label{cor:flow}
The unnormalized homogeneous Ricci flow starting from any
metric in Proposition~\ref{prop:scalar-construction} has finite
maximal forward time \(T\), and its scalar curvature tends to
\(+\infty\) as \(t\) tends to \(T\).
\end{corollary}
\begin{proof}
Apply \cite[Theorem~A]{AraujoFlow} to the compact-type base
\(\RR Z\oplus\mathfrak u\) and the abelian ideal \(V\).
There is one real-part weight, \(\alpha(aZ+X)=ac\); after
subtracting it the action is the semisimple purely imaginary
operator \(\pi(X)\). The metric makes the base and \(V\)
orthogonal, as required by the theorem. The simply connected
group is diffeomorphic to the product of the nontrivial
simply connected compact group with Lie algebra
\(\mathfrak u\) and \(\RR^{d+1}\). Its top compact-factor
cohomology is nonzero, so it is not Euclidean. Araujo's
theorem gives finite forward time. The scalar limit then
follows from \cite[Theorem~1.1(i)]{Lafuente}.
\end{proof}
\section{Sharp dimensions for a prescribed Levi algebra}
\label{sec:minimum-dimensions}

For a compact semisimple Lie algebra \(\mathfrak u\ne0\), let
\(\mu_{\mathbb R}(\mathfrak u)\) denote the least dimension of a
faithful real \(\mathfrak u\)-module. We fix a Levi algebra of the form
\[
  \mathfrak b=\mathfrak u\oplus\mathfrak l_-
    \oplus\mathfrak{sl}(2,\mathbb R)^{\oplus k},
  \qquad k\ge 0,
\]
where \(\mathfrak u\) is compact semisimple and \(\mathfrak l_-\) is
zero or is semisimple and admits a negative Ricci inner product.
Put \(\mathfrak b_{\mathrm{nc}}=\mathfrak l_-\oplus
\mathfrak{sl}(2,\mathbb R)^{\oplus k}\). The following minimization
allows arbitrary solvable radicals. At equality the positive-stable
direction is a consequence of negative Ricci curvature.

\begin{theorem}\label{thm:minimum-dimensions}
Among Lie algebras with Levi algebra \(\mathfrak b\) admitting a
negative Ricci inner product, the following dimensions and equality
structures hold.
\begin{enumerate}
\item If \(\mathfrak u\ne0\), the minimum dimension is
\[
  \dim\mathfrak b+1+\mu_{\mathbb R}(\mathfrak u).
\]
The Lie algebras attaining it are precisely those isomorphic to
\[
  \mathfrak b_{\mathrm{nc}}
  \oplus \bigl((\mathfrak u\oplus\mathbb R Z)\ltimes V\bigr),
  \qquad [Z,v]=Dv,
\]
where \(V\) is abelian and is a faithful real \(\mathfrak u\)-module of
dimension \(\mu_{\mathbb R}(\mathfrak u)\), and
\(D\in\operatorname{End}_{\mathfrak u}(V)\) has all its complex
eigenvalues in the open right half-plane. The line \(\mathbb R Z\) may
be oriented to arrange this sign. In each of these models \(V\) is the
actual nilradical.
\item If \(\mathfrak u=0\) and the radical is required to be nonzero,
the minimum dimension is \(\dim\mathfrak b+2\). The unique Lie algebra
up to isomorphism attaining it is
\(\mathfrak b\oplus\mathfrak{aff}(\mathbb R)\).
\end{enumerate}
\end{theorem}

\begin{proof}
Let \(\mathfrak r\) and \(\mathfrak n\) denote the radical and
nilradical. First, any unimodular Lie algebra with nonzero solvable
radical has a nonzero abelian ideal \(W\), obtained by taking the last
nonzero term of the derived series of \(\mathfrak r\). This term is an
ideal of the full algebra because the adjoint action preserves the
derived series. The mean curvature vector is zero by unimodularity.
Lemma~\ref{lem:abelian-trace} gives a nonnegative Ricci trace on \(W\)
for every inner product, which excludes negative Ricci curvature. In
particular, a nonzero radical in the theorem strictly contains its
nilradical. Indeed, if \(\mathfrak r=\mathfrak n\),
\(\operatorname{ad}_X\) is traceless for \(X\in\mathfrak r\) because
it maps the whole Lie algebra into \(\mathfrak r\) and acts nilpotently
there. The adjoint trace vanishes also on the perfect algebra
\(\mathfrak b\). This would give unimodularity.

Suppose \(\mathfrak u\ne0\). The kernel of the \(\mathfrak u\)-action on
\(\mathfrak n\) is a compact semisimple ideal of \(\mathfrak b\).
By Lemma~\ref{lem:radical}, \([\mathfrak b,\mathfrak r]\subseteq
\mathfrak n\). If \(\mathfrak k\) is the above kernel, then
\(\operatorname{ad}_X\operatorname{ad}_Y|_{\mathfrak r}=0\) for all
\(X,Y\in\mathfrak k\), so the action of the perfect algebra
\(\mathfrak k\) on \(\mathfrak r\) vanishes. Consequently
\(\mathfrak k\) is an ideal of the whole algebra. The compact ideal
obstruction (Lemma~\ref{lem:compact-obstruction}, cf.\
\cite[Lemma~1.6]{JP} and \cite[Proposition~2.1]{AraujoObstruction})
forces \(\mathfrak k=0\). Hence \(\dim\mathfrak n\ge
\mu_{\mathbb R}(\mathfrak u)\), and \(\mathfrak r\ne\mathfrak n\)
gives the stated lower bound. To attain it, choose a smallest faithful
real module \(V\), make it abelian, and take \(D=I_V\). The compact
extension has negative Ricci curvature by
Proposition~\ref{prop:scalar-construction}. Taking the orthogonal
direct sum with a negative Ricci metric on \(\mathfrak l_-\) and
applying Theorem~\ref{thm:attachment} \(k\) times yields a negative
Ricci inner product on the full algebra.

Now assume equality in this bound. Then \(\dim\mathfrak n=
\mu_{\mathbb R}(\mathfrak u)\) and \(\dim(\mathfrak r/\mathfrak n)=1\).
The induced \(\mathfrak u\)-action on
\(\mathfrak n/[\mathfrak n,\mathfrak n]\) is faithful. In fact, a
derivation acting trivially on that quotient also acts trivially on
each successive quotient of the lower central series, since these
quotients are spanned by the corresponding iterated brackets of
first-layer elements. Its action on \(\mathfrak n\) is therefore
nilpotent. For an element of \(\mathfrak u\), the same action is
skew-adjoint for a \(\mathfrak u\)-invariant inner product and so
vanishes if nilpotent. Faithfulness on \(\mathfrak n\) proves the
assertion. The minimal dimension of \(\mathfrak n\) now forces
\([\mathfrak n,\mathfrak n]=0\).

A minimal faithful real module has no trivial summand and no two
isomorphic irreducible summands, since deleting any such summand
preserves faithfulness. Thus for
\(\mathfrak n=\bigoplus_i W_i\) with \(W_i\) real irreducible,
\[
  \operatorname{End}_{\mathfrak u}(\mathfrak n)
    = \prod_i \operatorname{End}_{\mathfrak u}(W_i),
  \qquad
  \operatorname{End}_{\mathfrak u}(W_i)\simeq
    \mathbb R,\ \mathbb C,\ \text{or }\mathbb H.
\]
As a Lie algebra with the commutator bracket, the derived algebra of
this product is a direct sum of compact algebras
\(\mathfrak{sp}(1)\). The action of \(\mathfrak b_{\mathrm{nc}}\)
takes values in the commutant of \(\mathfrak u\) and hence, by
perfectness, in this compact derived algebra. Each of its simple
ideals is noncompact. For \(\mathfrak l_-\) this follows from the
compact ideal obstruction applied to its own negative Ricci metric.
A homomorphism from a noncompact simple real Lie algebra to a compact
Lie algebra is zero, since a nonzero such map would pull back an
invariant positive definite inner product. Therefore
\([\mathfrak b_{\mathrm{nc}},\mathfrak n]=0\).

The Levi action on \(\mathfrak r/\mathfrak n\) is trivial, and complete
reducibility supplies a \(\mathfrak b\)-invariant complementary line
\(\mathbb R Z\) in \(\mathfrak r\). This line centralizes
\(\mathfrak b\), so \(\mathfrak b_{\mathrm{nc}}\) is a direct factor
and \(D=\operatorname{ad}_Z|_{\mathfrak n}\) commutes with
\(\mathfrak u\). Each \(W_i\) is preserved by \(D\) and is an abelian
ideal of the full Lie algebra.

Consider any negative Ricci inner product and write its mean curvature
vector in the algebraic decomposition as
\(H=H_{\mathfrak b}+a Z+H_{\mathfrak n}\).
Lemma~\ref{lem:abelian-trace}, applied separately to each \(W_i\),
gives
\[
  0<\operatorname{tr}(\operatorname{ad}_H|_{W_i})
    =a\,\operatorname{tr}(D|_{W_i}).
\]
The Levi term has zero trace because its action on \(W_i\) is a
representation of a perfect algebra, and the nilradical term is zero
because \(\mathfrak n\) is abelian. After orienting \(Z\) so that
\(a>0\), each \(\operatorname{tr}(D|_{W_i})\) is positive. With respect
to a \(\mathfrak u\)-invariant reference inner product on \(W_i\),
the symmetric part of \(D|_{W_i}\) commutes with \(\mathfrak u\), so
it is scalar by irreducibility. Its scalar value is
\(\operatorname{tr}(D|_{W_i})/\dim W_i>0\). The skew-adjoint part
has purely imaginary spectrum, proving that every eigenvalue of
\(D\) has positive real part.

Conversely, in each displayed model \(V\) is the nilradical.
Any nonzero nilpotent ideal lies in the radical, and an element of the
radical with a nonzero \(Z\)-coefficient acts non-nilpotently on \(V\).
Such an element cannot belong to a nilpotent ideal. The model therefore
has the stated faithful action and the positive-stable direction
\(Z\) on its actual nilradical. Theorem~\ref{thm:mixed} gives negative
Ricci curvature and completes the equality characterization.

If \(\mathfrak u=0\) and \(\mathfrak r\ne0\), the same obstruction
to a nilpotent radical gives \(\dim\mathfrak r\ge2\). In dimension two
the nonnilpotent solvable Lie algebra is
\(\mathfrak{aff}(\mathbb R)\), with basis \(A,X\) and bracket
\([A,X]=X\). A derivation of this algebra maps \(A\) to a multiple
of \(X\) and maps \(X\) to a multiple of \(X\), as follows by applying
the derivation identity to \([A,X]=X\). Its derivation algebra is
solvable. Since \(\mathfrak b\) is perfect, its homomorphism to that
derivation algebra vanishes, giving the unique equality model. The
inner product making \(A,X\) orthonormal has Ricci operator \(-I\).
Combining it with the metric on \(\mathfrak l_-\) and using
Theorem~\ref{thm:attachment} gives existence at the bound.
\end{proof}

\begin{remark}\label{rem:minimum-dimensions-isomorphism}
For fixed \(\mathfrak b\), two equality models with \(\mathfrak u\ne0\),
described by \((V,\rho,D)\) and \((V',\rho',D')\), are isomorphic
if and only if there are \(\theta\in\operatorname{Aut}(\mathfrak u)\),
a real linear isomorphism \(T\colon V\to V'\), and \(c>0\) such that
\[
  T\rho(X)=\rho'(\theta X)T\quad(X\in\mathfrak u),
  \qquad TD=cD'T.
\]
Indeed, every isomorphism restricts to such a map \(T\) on the actual
nilradical and takes \(Z\) to \(cZ'+v'\) for some \(c\ne0\) and
\(v'\in V'\), because it preserves the radical as well. The induced
isomorphism on the quotient by the nilradical preserves the derived
algebra and its compact ideal, so it determines \(\theta\). Any
discrepancy between the image of a chosen Levi lift of \(X\) and the
chosen lift of \(\theta X\) belongs to \(V'\). Since \(V'\) is abelian,
neither this discrepancy nor \(v'\) contributes to the action on
\(V'\). The bracket then gives the two identities, and the common
positive-real-part convention for \(D,D'\) forces \(c>0\).
Conversely, these identities define an isomorphism by acting as
\(\theta\) on \(\mathfrak u\), as \(T\) on \(V\), and by sending
\(Z\) to \(cZ'\) and identifying the fixed noncompact factors.
\end{remark}

\begin{corollary}\label{cor:small-compact-dimensions}
For the following prescribed compact simple Levi algebras, the exact
minimum dimensions are as follows.
\[
\begin{array}{c|c|c|c|c}
\text{Levi algebra }\mathfrak u & \dim\mathfrak u
 & \text{minimal faithful real module}
 & \mu_{\mathbb R}(\mathfrak u) & \min\dim\mathfrak g\\ \hline
\mathfrak{su}(2)&3&\mathfrak{su}(2)\text{ (adjoint)}&3&7\\
\mathfrak{su}(3)&8&(\mathbb C^3)_{\mathbb R}&6&15\\
\mathfrak{so}(5)&10&\mathbb R^5&5&16\\
\mathfrak g_2\text{ (compact)}&14&\operatorname{Im}\mathbb O&7&22
\end{array}
\]
\end{corollary}

\begin{proof}
A faithful real \(d\)-dimensional module of a compact semisimple
Lie algebra gives an embedding in \(\mathfrak{so}(d)\) by choosing an
invariant inner product. Dimensions show that the adjoint module of
\(\mathfrak{su}(2)\) and the vector module of \(\mathfrak{so}(5)\) are
minimal. For \(\mathfrak{su}(3)\), Weyl's dimension formula
\(
  \dim V_{a,b}=(a+1)(b+1)(a+b+2)/2
\)
shows that the only nontrivial complex irreducibles of dimension at
most five are the standard module and its dual. These are distinct
complex conjugate modules. In the complexification of any real
module they must occur with equal multiplicities, so a nontrivial
real module has dimension at least six. The realification of the
standard module attains six.

The compact \(\mathfrak g_2\) acts faithfully on the
seven-dimensional space of imaginary octonions \cite[Theorem~4]{Baez}.
Dimension excludes an embedding in \(\mathfrak{so}(d)\) for \(d\le5\).
An embedding in \(\mathfrak{so}(6)\) would give a codimension-one
subalgebra. For an invariant positive definite inner product on
\(\mathfrak{so}(6)\), the orthogonal line would be preserved by that
subalgebra and hence centralized by it, since its action on the line
is skew-adjoint. The line would then be central in
\(\mathfrak{so}(6)\), which is impossible. Applying
Theorem~\ref{thm:minimum-dimensions} yields the table.
\end{proof}

\begin{remark}
The seven-dimensional \(\mathfrak{su}(2)\) example is due to Will
and is recorded in \cite[Section~2.2]{LWSurvey}.
The same survey's Table~1 records examples built from complex
representations regarded as real modules. The real vector module of
\(\mathfrak{so}(5)\) and the real octonionic module of the compact
\(\mathfrak g_2\) each use one real form of the corresponding complex
module. The dimension bounds in
Theorem~\ref{thm:minimum-dimensions} hold with arbitrary radicals.
\end{remark}
\section{Low-dimensional classification and explicit metrics}\label{sec:lowdim}

\subsection{The five-dimensional contact metric family}

Boeckx's contact metric models for \(I<-1\) in dimension five
are locally based on
\(\mathfrak{sl}(2,\mathbb R)\oplus\mathfrak{aff}(\mathbb R)\),
as observed in \cite[Remark~1]{LL}. Loiudice and Lotta
construct a separate family on this algebra
\cite[Theorem~6 and Corollary~1]{LL}. The Ricci formula
and a positive \(\mathcal D\)-homothety give a negative
Ricci metric in each fixed class \(I<-1\) (see below).
We compute the exact negative Ricci interval within the
unrescaled family of Loiudice and Lotta.
In the basis \(e_1,\ldots,e_5\) of their equation (6), the change of basis
\[
 K=-2e_5,\quad P=\sqrt2e_3,\quad H=\sqrt2e_4,\qquad
 A=e_1/\sqrt2,\quad N=e_2-\sqrt2e_5
\]
gives \([K,P]=-2H\), \([P,H]=2K\), \([H,K]=-2P\), and \([A,N]=N\).
The two spans commute. Thus the original bracket has the indicated direct
sum, and its nilradical is \(\mathbb RN\).

Use the positive parameter \(a>1\) of their theorem statement. The
parameter in the proof's equation (7) is \(-a\). Set
\(\lambda=(a-a^{-1})/4\) and \(F=\sqrt{1+4\lambda^2}\). The contact metric
parameters are \(\kappa=1-\lambda^2\) and \(\mu=2+F\). In dimension five
the Ricci operator formula for non-Sasakian contact \((\kappa,\mu)\) metrics
gives the eigenvalues
\[
 4(1-\lambda^2),\qquad
 2-2\mu+(2+\mu)\lambda,\qquad
 2-2\mu-(2+\mu)\lambda,
\]
with multiplicities \(1,2,2\), respectively
\cite[Equations (2.6)--(2.7)]{LottaMM}. The Reeb eigenvalue requires
\(\lambda>1\). The larger horizontal eigenvalue is
\(f(\lambda)=4\lambda-2+(\lambda-2)F\). For \(1\leq\lambda\leq2\),
\(4\lambda/F<2\) gives \(f'(\lambda)>0\), and
\(f(1)=2-\sqrt5<0<f(2)\). For \(\lambda\geq2\), the formula for \(f\)
gives \(f(\lambda)>0\) directly. On \((1,2)\) both sides of
\((2-\lambda)F=4\lambda-2\) are positive, and squaring reduces this equation
to \(4\lambda^3-16\lambda^2+\lambda+12=0\). Thus the unrescaled family has
negative Ricci curvature exactly when \(1<\lambda<\lambda_*\), where
\(\lambda_*\) is the unique root in \((1,2)\) of this cubic.

More generally, for fixed Boeckx invariant \(I=-b<-1\), a positive
\(\mathcal D\)-homothety can prescribe any new value \(\lambda'>0\).
One then has \(\mu'=2+2b\lambda'\), so the new metric has negative Ricci
curvature exactly when
\[
 1<\lambda'<\frac{b-1+\sqrt{(b-1)^2+b}}{b}.
\]
Indeed, the upper horizontal eigenvalue is
\(2[b(\lambda')^2+2(1-b)\lambda'-1]\), which is negative at \(\lambda'=1\)
for \(b>1\). For instance, \(a=2\) followed by the \(\mathcal D\)-homothety
with factor \(3/10\) gives \(\lambda'=5/4\) and \(\mu'=37/6\), and hence
the three Ricci eigenvalues \(-9/4,-1/8,-493/24\).

\subsection{Classification in dimensions five and six}

\begin{theorem}\label{thm:lowdim}
Among nonsolvable, nonsemisimple real Lie algebras admitting a strictly
negative Ricci inner product, the smallest dimension is five, with unique
isomorphism type
\[
 \mathfrak{sl}(2,\mathbb R)\oplus\mathfrak{aff}(\mathbb R).
\]
In dimension six the isomorphism types are precisely the standard affine
algebra
\(\mathfrak{gl}(2,\mathbb R)\ltimes\mathbb R^2\) and the direct sums
\[
 \mathfrak{sl}(2,\mathbb R)\oplus
    \bigl(\mathbb RA\ltimes_D\mathbb R^2\bigr),
 \qquad [A,v]=Dv,\qquad
 \operatorname{Re}\operatorname{Spec}D\subset(0,\infty).
\]
In the latter family, two matrices give the same isomorphism type precisely
when they are similar up to multiplication by a positive real scalar. With
\(\operatorname{tr}D=2\), representatives are
\[
 \begin{pmatrix}1+a&0\\0&1-a\end{pmatrix}\quad(0\leq a<1),
 \qquad
 \begin{pmatrix}1&1\\0&1\end{pmatrix},
 \qquad
 \begin{pmatrix}1&-b\\b&1\end{pmatrix}\quad(b>0).
\]
\end{theorem}

\begin{proof}
Let \(\mathfrak r\ne0\) be the solvable radical and \(\mathfrak s\ne0\)
a Levi factor. In dimension at most six, \(\dim\mathfrak s\leq5\).
The only complex semisimple Lie algebra of positive dimension at most five
is \(\mathfrak{sl}(2,\mathbb C)\). Indeed, rank one gives the root system
\(A_1\), and a complex semisimple Lie algebra of rank at least two has
dimension at least six. Thus \(\mathfrak s\) is either
\(\mathfrak{so}(3)\) or \(\mathfrak{sl}(2,\mathbb R)\), both of dimension
three.

If \(\mathfrak r\) is unimodular, then so is the full algebra:
\(\operatorname{ad}X\) for \(X\in\mathfrak r\) has zero trace on the
quotient by \(\mathfrak r\), and the trace character vanishes on the
perfect algebra \(\mathfrak s\). The final nonzero term of the derived
series of \(\mathfrak r\) is a nonzero abelian ideal of the full algebra.
For a unimodular algebra the mean-curvature vector vanishes, so
Lemma~\ref{lem:abelian-trace} excludes negative Ricci curvature.
In particular this excludes every radical of dimension at most one.

In dimensions five and six the remaining radicals have dimensions two
and three. Let \(\tau(X)=\operatorname{tr}(\operatorname{ad}_{\mathfrak r}X)\).
Since \(\tau\ne0\), its kernel \(V\) is an ideal of codimension one.
It is abelian when \(\dim V=1\). When \(\dim V=2\), take a basis \(v,w\)
and write \([v,w]=\alpha v+\beta w\). The traces on \(V\) of
\(\operatorname{ad}_{\mathfrak r}v\) and
\(\operatorname{ad}_{\mathfrak r}w\) are \(\beta\) and \(-\alpha\).
They agree with the traces on \(\mathfrak r\), because both adjoint
operators act trivially on \(\mathfrak r/V\). Thus \(\alpha=\beta=0\).
It follows that \(V\) is the actual nilradical of \(\mathfrak r\) and of
the full algebra: \(V\) is nilpotent, and the only larger subspace of
\(\mathfrak r\) containing \(V\) is the nonunimodular algebra \(\mathfrak r\)
itself.

The induced \(\mathfrak s\)-action on \(\mathfrak r/V\) is trivial by
perfectness. Complete reducibility allows a complement \(\mathbb RA\)
to \(V\) with \([\mathfrak s,A]=0\). Choose the sign of \(A\) so that
\(D=\operatorname{ad}_A|_V\) has positive trace. Write
\(\pi:\mathfrak s\to\mathfrak{gl}(V)\) for the remaining action.
The Jacobi identity gives \([D,\pi(X)]=0\) for \(X\in\mathfrak s\).
If \(\dim V=1\), then \(\pi=0\). A compact Levi factor acts trivially
also when \(\dim V=2\), because its representation is conjugate into
\(\mathfrak{so}(2)\) and \(\mathfrak s\) is perfect. A trivial action
of the compact factor makes it a compact semisimple ideal, which is
excluded by Lemma~\ref{lem:compact-obstruction}.

Consequently the Levi factor is split. When \(\dim V=1\),
\(\mathfrak r\cong\mathfrak{aff}(\mathbb R)\), and the five-dimensional
existence follows as above or from Theorem~\ref{thm:attachment}. When
\(\dim V=2\) and \(\pi=0\), Theorem~\ref{thm:attachment-reverse} gives exactly the
positive-real-part condition on \(D\): the adjoint action of any
\(xA+v\in\mathfrak r\) on \(V\) is \(xD\), and \(\operatorname{tr}D>0\)
forces \(x>0\) if all the real parts of the eigenvalues of \(xD\) are
positive.

In the remaining case \(\pi\ne0\). Its image is all of
\(\mathfrak{sl}(V)\), since it is injective by simplicity and its image
has trace zero and dimension three. Its centralizer in
\(\mathfrak{gl}(V)\) consists of scalars. After rescaling \(A\), we have
\(D=I\). Identifying \(\mathfrak s\) with its image gives the standard affine
algebra. Its negative Ricci metric is given below.

Finally, an isomorphism between two members of the direct-sum family
preserves the characteristic radical and nilradical. If it acts as
\(T\) on \(V\) and sends \(A\) to \(xA'+v'\), then
\(TD=xD'T\). Both traces are positive, hence \(x>0\). Conversely,
this identity defines an isomorphism on the radicals and hence on the
direct sums. The normal forms are the real Jordan forms after positive
trace normalization. The standard affine algebra is distinct from all
these direct sums: the image of the full adjoint action on \(V\) has
dimension four for the affine algebra and dimension one for a direct sum.
\end{proof}

\subsection{Explicit metrics in dimensions six and seven}

\begin{example}[The six-dimensional affine algebra]\label{ex:active-affine}
Write the affine algebra as
\[
 \mathfrak g=(\mathbb R A\oplus\mathfrak{sl}(2,\mathbb R))\ltimes V,
 \qquad V=\mathbb R^2.
\]
The simple algebra acts by its standard representation and \([A,v]=v\).
In the simple factor take the basis
\[
 K=\begin{pmatrix}0&-1\\1&0\end{pmatrix},\quad
 P=\begin{pmatrix}0&1\\1&0\end{pmatrix},\quad
 H=\begin{pmatrix}1&0\\0&-1\end{pmatrix}.
\]
Take \(K,P,H,4A\) orthonormal and perpendicular to \(V\), and give \(V\)
metric matrix \(\operatorname{diag}(4,1/4)\) in its original standard
basis. The action matrices in the resulting orthonormal basis of \(V\)
are obtained by conjugating by \(h=\operatorname{diag}(2,1/2)\).
They are
\[
 A_K=\begin{pmatrix}0&-4\\ 1/4&0\end{pmatrix},\qquad
 A_P=\begin{pmatrix}0&4\\ 1/4&0\end{pmatrix},\qquad A_H=H.
\]
The intrinsic Ricci eigenvalues on \(K\), \(P\), and \(H\) are \(2\),
\(-6\), and \(-6\). The symmetric parts of \(A_K,A_P\) have off-diagonal
entries \(-15/8,17/8\). Subtracting their Gram matrix gives the
\((K,P)\) Ricci block
\[
 \frac1{32}\begin{pmatrix}-161&255\\255&-481\end{pmatrix},
 \qquad \det=\frac{97}{8}.
\]
The \(H\) and \(4A\) entries are \(-8\) and \(-32\), respectively. On
\(V\), the commutator and mean-curvature terms give
\[
 \frac12\sum_{X=K,P,H}[A_X,A_X^*]-32I
 =\frac1{16}\operatorname{diag}(-257,-767).
\]
The mean-curvature vector is \(8(4A)\). All remaining cross terms
vanish: the involution acting by \(-I\) on \(V\) and by \(I\) on its
complement is isometric, all simple-factor actions are trace free, and
the symmetric \(H\)-action is Frobenius-orthogonal to the two other
symmetric parts. The displayed block and the remaining entries are
strictly negative.
\end{example}

\begin{example}[Trivial simple-factor action in dimensions six and seven]
\label{ex:inert-scalar-heisenberg}
First let \(\mathfrak r_d=\mathbb RA\ltimes\mathbb R^d\) with
\([A,x]=x\). In its orthonormal basis \(A,X_1,\ldots,X_d\), the
Ricci operator is \(-dI\). Thus Theorem~\ref{thm:attachment} gives
examples of dimensions \(3m+d+1\) for all \(d,m\geq1\).
For \(d=2,m=1\), in that theorem's frame \(k,P,q\) with \(P=p+tv\),
choose \(v=X_1\) and \(t=1/10\). The \((P,X_1)\) block is
\[
 \begin{pmatrix}
 -9/160&-1/5\\-1/5&-199/100
 \end{pmatrix},\qquad
 \det=\frac{1151}{16000}.
\]
The remaining block in the order \(k,q,A,X_2\) is
\[
 \operatorname{diag}\bigl(-3/800,-9/160,-401/200,-2\bigr).
\]
To check the block substitutions, the
mean-curvature vector of \(\mathfrak r_2\) is \(2A\), and
\(L=\operatorname{ad}_{\mathfrak r_2}X_1\) sends \(A\) to \(-X_1\)
and vanishes on \(\mathbb R^2\).

For a nonabelian nilradical, take
\(\mathfrak r=\mathbb RA\ltimes\mathfrak h_3\) with orthonormal basis
\(A,X,Y,Z\) and brackets
\[
 [X,Y]=Z,\qquad [A,X]=X,\qquad [A,Y]=Y,\qquad [A,Z]=2Z.
\]
The mean-curvature vector is \(4A\). The nilpotent Ricci operator and
the mean-curvature term give the block on \(X,Y,Z\),
\[
 \tfrac12\operatorname{diag}(-1,-1,1)-4\operatorname{diag}(1,1,2).
\]
The \(A\)-direction contributes \(-6\), and its mixed terms
vanish since \(\operatorname{ad}_A|_{\mathfrak h_3}\) is diagonal and the
operators \(\operatorname{ad}_w|_{\mathfrak h_3}\),
\(w\in\mathfrak h_3\), have zero diagonal in this basis. Thus, in the order
\(A,X,Y,Z\), we have
\[
 \operatorname{Ric}_{\mathfrak r}=\operatorname{diag}(-6,-9/2,-9/2,-15/2).
\]
Apply Theorem~\ref{thm:attachment} with \(v=Z\) and \(t=1/10\).
The operator \(L=\operatorname{ad}_{\mathfrak r} Z\) sends \(A\) to
\(-2Z\) and vanishes on \(\mathfrak h_3\). The \((P,Z)\) Ricci block is
\[
 \begin{pmatrix}
 -57/800&-13/20\\-13/20&-299/40
 \end{pmatrix},\qquad
 \det=\frac{3523}{32000}.
\]
The remaining block, in the order \(k,q,A,X,Y\), is
\[
 \operatorname{diag}\bigl(-3/800,-9/160,-301/50,-9/2,-9/2\bigr).
\]
Both metric Lie algebras
therefore have strictly negative Ricci curvature.
\end{example}
\section{Further questions}\label{sec:further}

Theorem~\ref{thm:attachment-reverse} leaves open whether
existence on \(\mathfrak{sl}(2,\RR)^{\oplus m}\oplus\mathfrak r\)
always implies existence on the solvable factor
\(\mathfrak r\). The trace condition in
Lemma~\ref{lem:abelian-trace} applies to any abelian ideal
invariant under the adjoint action. The three classes
considered here have enough such ideals to recover their
solvable criteria. For other nilradicals one may ask what
additional restrictions transfer from the full direct sum
to its solvable factor.

Theorem~\ref{thm:mixed} concerns stable radical directions
and a specified class of Levi algebras. A general existence
criterion for noncompact Levi factors remains a separate
problem. The eleven-dimensional example also leaves the
minimum dimension of a semisimple counterexample to
factorwise inheritance undetermined.

\end{document}